\documentclass[11pt,a4paper]{amsart}

\usepackage{amscd,amssymb,amsthm}

\title[Existence of indecomposable rank two vector bundles]{Existence of indecomposable rank two vector bundles on higher dimensional toric varieties}
\author{Giulio Cotignoli and Alexandru Sterian}
\date{}

\usepackage{amsfonts, amssymb, amscd}
\usepackage{verbatim}
\usepackage{eucal}
\usepackage{amssymb}
\usepackage{mathrsfs}

\def\OO{{\mathcal O}}
\def\LL{{\mathcal L}}

\def\ZZ{{\mathbb Z}}

\addtocounter{footnote}{1}

\newtheorem{lemma}{Lemma}[section]
\newtheorem{theorem}[lemma]{Theorem}

\newtheorem{conjecture}[lemma]{Conjecture}
\newtheorem{proposition}[lemma]{Proposition}
\theoremstyle{definition}
\newtheorem{definition}[lemma]{Definition}

\newtheorem{remark}[lemma]{Remark}

\theoremstyle{remark}
\newtheorem*{proof*}{Proof}

\begin{document}
\maketitle

\begin{abstract}
In the mid 70's, Hartshorne conjectured that, for all $n > 7$, any rank $2$ vector bundles on $\mathbb{P}^n$ is a direct sum of line bundles. This conjecture remains still open. In this paper, we construct indecomposable rank two vector bundles on a large class of Fano toric varieties. Unfortunately, this class does not contain $\mathbb{P}^n$.
\end{abstract}

 \bigskip

\centerline{\textbf{AMS classification [2000]} : {14M25, 14J60}}

\section{Introduction}

The starting point of our problem is Hartshorne conjecture
concerning the existence of indecomposable rank two vector bundles
on projective spaces \cite{Ha}:

\begin{conjecture}Every rank $2$ vector bundle on $\mathbb{P}^n$, $n>7$ is a
direct sum of line bundles.
\end{conjecture}
On $\mathbb{P}^3$ there are plenty of indecomposable rank $2$
vector bundles. They have been studied with the hope that a good
understanding of the situation for $n = 3$ would have been for
some help for $n>7$. Up to now, no indecomposable rank $2$ vector
bundles have been found on $\mathbb{P}^n, n>4$. On $\mathbb{P}^4$,
we know only one type of indecomposable rank $2$ vector bundles,
the Horrocks-Mumford bundles. More recently, people have tried to
find indecomposable rank $2$ vector bundles on any smooth variety.
For instance, Laura Costa and Rosa Maria Mir\'o-Roig prove in \cite{Co-MR} the existence
of rank two vector bundles on projective bundles over algebraic
curves.

In this paper, using a similar techniques as in \cite{Co-MR}, we construct
indecomposable rank $2$ vector bundles on large class of Fano toric varities.

\bigskip

The main theorem of this paper is the following:

\begin{theorem} Let $X$ be a Fano pseudo-symmetric  complete toric variety of dimension $n>2$, other than a product of ${\mathbb{P}^1}$. There exists an indecomposable rank $2$ vector bundle on $X$.
\end{theorem}

Note that this result does not apply to $\mathbb{P}^n$.

\section{About pseudo-symmetric toric varieties.}
Let $X$ be a smooth, complete toric variety of dimension $n$. Let
$T$ be its complex torus, $N=Hom_{\mathbb{Z}}(\mathbb{C}^*,T)$,
$N_{\mathbb{Q}}=N\otimes_{\mathbb{Z}}\mathbb{Q}$ and $\Sigma_X$
the fan of $X$ in $N_{\mathbb{Q}}$. We denote by $\langle x_1,...,x_h \rangle$
the convex cone spanned by $x_1,...,x_h\in N$ in $N_{\mathbb{Q}}$.

\begin{definition}
Let $X$ be complete toric variety, then the set of ray generators
of $X$ is:
$$G(\Sigma_X)=\{x\in N| \langle x \rangle \in\Sigma_X \,\,\textrm{and}\,\, x \,\,\textrm{is primitive in}\,\, \langle x \rangle \cap N\}.$$
\end{definition}

\begin{definition}
A toric Fano variety is called pseudo-symmetric if its
fan containts two centrally symmetric maximals cones, i.e. there exists $\sigma, \sigma'\in\Sigma_X$ maximals cones such that $\sigma=-\sigma'$.
\end{definition}

\bigskip

Let $e_1,...,e_n$ be a basis of the lattice $\mathbb{Z}^n$ and suppose
that $n=2r$ is even. Consider the following elements of
$\mathbb{Z}^n$:

\bigskip

\noindent
\begin{tabular}{l}
$v_i = e_i,\, \textrm{for} \,\, i=1,...,n$, \\
$v_{n+1}=-e_1-...-e_n$,\\
$v_{i+n+1}=-v_i\,\,\textrm{for}\,\,i=1,..,n+1.$\\
\end{tabular}

\begin{definition}
The $n$-dimensional \emph{Del Pezzo variety} $V^n$ is the toric variety associated to the fan $\Sigma_{V^n}\in N_{\mathbb{Q}}.$
such that
$$G(\Sigma_{V^n})=\{v_1,...,v_{2n+2}\}.$$
In the same way, the $n$ dimensional \emph{pseudo Del Pezzo variety}
$\widetilde{V^n}$ is the toric Fano variety described by the fan
$\Sigma_{\widetilde{V}^n}\subset N_{\mathbb{Q}}$ such that
$$G(\Sigma_{\widetilde{V}^n})=\{v_1,...,v_{2n+1}\}.$$
\end{definition}

\begin{remark}
The Picard number for the $n$-dimensional Del Pezzo variety $V^n$ is $\rho_{V^n}=n+2$. For $n=2$, this variety is the Del
Pezzo surface $S_3$, obtained by blowing up
$\mathbb{P}^2$ at three points.
\end{remark}

Del Pezzo varieties play a very important role in the description
of Fano toric varieties. Actually, Del
Pezzo varieties were introduced by V.E. Voskresenskii and A.A.
Klyachko, and they show that every toric Fano variety whose fan is
centrally symmetric is a product  of Del Pezzo varieties and a
product of $\mathbb{P}^1$. This result has been generalized by G.
Ewald in the next theorem:

\begin{theorem}( see \cite{Ew})
For any pseudo-symmetric Fano variety $X$, exists $s, m, n\in
\mathbb{Z}_{\geq 0}$ and $k_1,...,k_m, l_1,...,l_n \in\mathbb{Z}_{\geq 0}$
such that:
$$X\simeq(\mathbb{P}^1)^s\times V^{2k_1}\times...\times V^{2k_m}\times(\widetilde{V})^{2l_1}\times...\times(\widetilde{V})^{2l_n}.$$
\end{theorem}

\bigskip
\begin{remark} Let $e_1,...,e_n$ be a basis of $\mathbb{Z}^n$ and let $\Sigma_{\mathbb{P}^n}$ be the fan of $\mathbb{P}^n$ in $\mathbb{Z}^n$. Then, we have:
$$ G(\Sigma_{\mathbb{P}^n}) = \{e_1,...,e_n,-e_1-...-e_n \}.$$
It's obvious, that for $n \geq 2$, $\mathbb{P}^n$ does not contain two centrally symmetric maximals cones, so unfortunately the main theorem does not apply to $\mathbb{P}^n$.
\end{remark}
In the last section of this paper, we will prove the existence of an indecomosable rank $2$ vector bundle on any Del Pezzo variety of dimension $n>2$.

\section{Cohomology of Line Bundles on Del-Pezzo Varieties}
We follow here the approach of Borisov and Hua \cite{BoH}. They
provide a description of cohomology of a line bundle $\mathcal{L}$
on every toric variety. We will use their construction for
Del-Pezzo varieties.

\bigskip
Let $V^n$ be the Del-Pezzo toric variety of dimension $n$ ($n$
even). For every $\textbf{r}=(r_i)_{i=1,...,2n+2} \in \mathbb{Z}^{2n+2}$,
we denote by $Supp(\textbf{r})$ the simplicial complex on $2n+2$
vertices $\{1,..,2n+2\}$ which consists of all subsets
$J\subseteq\{1,..,2n+2\}$ such that $r_i\geq 0$ for all $i\in J$
and there exists a cone in $\Sigma_{V^n}$ that contains all the
ray generators $v_i$, $i\in J$. We will abuse notation and also
denote by $Supp(\textbf{r})$ the subfan of $\Sigma_{V^n}$ whose
cones are the minimal cones of $\Sigma_{V^n}$ that contain all
$v_i,i\in J$ for all subsets $J$ as above. It should be clear from
the context whether $Supp(\textbf{r})$ refers to the simplicial
complex or to its geometric realization as a subfan of
$\Sigma_{V^n}$.

For example, if all coordinates $r_i$ are negative then the
simplicial complex $Supp(\textbf{r})$ consists of the empty set
only, and its geometric realization is the zero cone of
$\Sigma_{V^n}$. In the other extreme case, if all $r_i$ are
nonnegative then the simplicial complex $Supp(\textbf{r})$ encodes
the fan $\Sigma_{V^n}$, which is its geometric realization.

\begin{proposition}\label{omologia}
Let $V^n$ be the Del-Pezzo toric variety of dimension $n$ and let
$E_i$ be the toric divisors associated to the ray generators $v_i$
for $i=1,...,2n+2$. The cohomology $H^p(V^n,\mathcal{L})$ is
isomorphic to the direct sum over all
$\mathbf{r}=(r^i)_{i=1,...,2n+2}$, such that
$\mathcal{O}(\sum_{i=1}^{2n+2} r_iE_i)\cong\mathcal{L}$, of the
($n-p$)-th reduced homology of the simplicial complex
Supp($\textbf{r}$).
\end{proposition}

See \cite{BoH} for a proof.

\begin{remark}
For example, $H^0(V^n,\mathcal{L})$ only comes from $\bf r$ for
which ${\rm supp}({\bf r})$ is the entire fan $\Sigma$, i.e.
$\OO(\sum_{i=1}^{2n+2} r_i E_i)\cong \LL$ with ${\bf r}\in
\ZZ_{\geq 0}^{2n+2}$. In the other extreme case $
H^n(V^n,\mathcal{L})$ only appears when the simplicial complex
${\rm supp}({\bf r})=\{\emptyset\}$, i.e. when
$\OO(\sum_{i=1}^{2n+2} r_i E_i)\cong \LL$ with all $r_i\leq -1$.
\end{remark}
\begin{remark}
If $Z=X\times Y$ is the product variety of $X$ and $Y$, $E$ an indecomposable rank two vector bundle on $Y$ and $p:Z\longrightarrow Y$ the projection on $Y$, 
then the pull-back  $p^*(E)$  is still an indecomposable rank 2 vector bundle.
\end{remark}
As a consequence of the above remark, the existence of an indecomposable rank $2$ vector bundle on any Del Pezzo and pseudo Del Pezzo variety of dimension $n>2$, will imply the existence of an indecomposable rank $2$ vector bundle on any pseudo-symmetric toric Fano variety $X$, other then a product of  ${\mathbb{P}^1}$.
\section{Indecomposable rank $2$ vector bundles on Del Pezzo
varieties} Our goal is the proof of the following fact:
\begin{theorem}
Let $V^n$ be the $n$-dimensional Del Pezzo toric variety. For every $n>2$, there exists an
indecomposable rank $2$ vector bundle on $V^n$.
\end{theorem}
An analogous proof show the existence of an indecomposable rank-2 vector bundle on pseudo-Del Pezzo variety.

\bigskip

 Given two line bundles $\mathcal{L}_1$,
$\mathcal{L}_2$ on $V^n$, any extension $\mathcal{E}$ of
$\mathcal{L}_1$ by $\mathcal{L}_2$  is a rank two vector bundle on
$V^n$ such that there is an exact sequence:
$$0\longrightarrow\mathcal{L}_1\longrightarrow\mathcal{E}\longrightarrow\mathcal{L}_2\longrightarrow 0.$$
Let $\mathcal{L}\sim{\mathcal{L}_2}^\nu\otimes\mathcal{L}_1$. It
is a standard fact that any extension of $\mathcal{L}_1$ by
$\mathcal{L}_2$ appears as a class in $H^1(V^n,\mathcal{L})$. As a
consequence, $\mathcal{L}_1$ and $\mathcal{L}_2$  give a non
trivial extension if and only if $H^1(V^n,\mathcal{L}) \neq 0$.
Thus, we want a vector bundle $\mathcal{E}$ whose class in
$H^1(V^n,{\mathcal{L}_2}^\nu \otimes \mathcal{L}_1$) is not
trivial for $\mathcal{L}_1$ and $\mathcal{L}_2$, such that
$\mathcal{E}$ is an extension of $\mathcal{L}_1$ by
$\mathcal{L}_2$.

\bigskip

First, we want to find a line bundle $\mathcal{L}$ on $V^n$ with
$H^1(V^n,\mathcal{L})\neq 0$. Let $E_1,...,E_{2n+2}$ be the
divisors corresponding to the ray generators $v_1,...,v_{2n+2}$ in
the fan of $V^n$. Recall that $\textrm{Pic}(V^n)$ is generated by
$E_1,...,E_{2n+2}$, so every line bundle $\mathcal{L}$ on $V^n$
can be written:
$$\mathcal{L}\sim \mathcal{O}(\sum_{i=1,...,2n+2} r_i E_i),$$ for some integers $r_1,...,r_{2n+2}$.

Now Proposition $\ref{omologia}$ lead us to restate our problem.
We look for an divisor $D=\sum_{i=1,...,2n+2} r_i E_i$  such such that the corresponding simplex
$Supp(\textbf{r})$ has $H_{n-1}(\emph{Supp}(\textbf{r}), \mathbb{C})\neq 0$.
Elements in this groups are $(n-1)$-cycles in
$\emph{Supp}(\textbf{r})$, that is, simplices which can be
described in a purely combinatorial way from the fan of $V^n$.
Using Proposition $\ref{omologia}$, we consider the simplicial complex $C_{\Sigma}(\textbf{r})$ associated
to the fan $\Sigma$ corresponding to the Del Pezzo variety $V^n$.
\newline We take $C_{-1}=\mathbb{C}$, since $C_{-1}$ encodes the zero cone.
\newline For every $i\in\{0,...,n\}$ we put $C_i$ the $\mathbb{C}$-vector space of all the $i$-simplices in the complex. 
\newline This simplices encodes the $(i+1)$-dimensional cones $<v_{j_1},...,v_{j_{i+1}}>$ which appear in the fan $\Sigma$, such that $r_{j_k} \geq 0$ for all $k \in \{1,...,i+1\}$.
\newline The differential operators $\delta_i:C_i\longrightarrow C_{i-1}$  are defined in the following way:
\newline $\delta_{-1}$ is the zero map;
\newline $\delta_0(<v_{j_k}>)=1, \forall j_k\in\{1,...,2n+2\}$; this means that every $0$-simplex is send to 1; the geometrical meaning of this aplication being that every $1$-dimensional cone is send in the zero cone;
\newline To describe $\delta_i$, with $i\in\{1,...,n\}$, we start by defining an order $j_1<...<j_{n}$. 
\newline Then, $\delta_i (<v_{j_1},...,v_{j_{i+1}}>)=\sum_{k=1}^n(-1)^k<v_{j_1},...,\hat{v}_{j_k},...,v_{j_{i+1}}>$, where the "hat" symbol over $v_{j_k}$ indicates that we take the $(i-1)$-simplex which doesn't contain the vertex $v_{j_k}$. So, this application is up to a sign, the restriction map that send an $i$-simplex in all its faces: $(i-1)$-simplices. Moreover, this map is equivalent with the one from $\check{C}$ech complex of $V^n$ Del Pezzo variety that send an $(i+1)$-dimensional cone in all its faces, $i$-dimensional cones that dont contain $v_{j_k}$, for $k\in\{1,...,i+1\}$ with the sign $(-1)^{k}$. For more details, see the proof of Proposition 4.1 from \cite{BoH}.

\bigskip

We are lead to study the fan $\Sigma_{V^n}$ of $V^n$. Starting
from the ray generators $v_1,...,v_{2n+2}$ in $\Sigma_{V^n}$,  split
them in two subsets:
$$x_i=v_i, \, \textrm{for} \, 1\leq i \leq n+1,$$
and
 $$y_j=v_{j+n+1}, \, \textrm{for} \, 0\leq j \leq n+1.$$
This splitting agrees with the structure of
our fan $\Sigma_{V^n}$. Indeed, the fan $\Sigma_{V^n}$ is defined
by prescribing that the rays $y_i$ and $x_i$ can not
appear at the same time in any cone that $\Sigma_{V^n}$ contains
as part of its support (see \cite{Ca}). More precisely, for $n=2r$, the fan $\Sigma_{V^n}$ is the union of the $\Sigma(m)$ for $0 \leq m \leq r$, where:
$$\Sigma (m)=\{\langle x_i,y_j \rangle_{i\in I, j\in J} | I, J\subset \{1,2,...,n+1\}, I \cap J = \emptyset,$$
$$\#I \leq r, \#J \leq r, \#I\cup J=m\}.$$

In our case, it be more simple to look for a divisor $D$ such that the first group of reduced homology $H_1(D)\neq 0.$
Using Proposition \ref{omologia}, this will be equivalent with $H^{n-1}(V^n,D)\neq 0,$ and so applying Serre duality we have $H^{1}(V^n,K_X-D)^{\vee}\neq 0.$

Let's consider the divisor $D=K_{V^n}+2t(E_p+E_{2n+2})$, with $t\in\mathbb{Z}$, $t\geq 1$, $p$ a fixed integer such that, $1\leq p\leq n$ and $K_{V^n}=-\sum_{k=1}^{2n+2}E_i$ the canonical divisor of $V^n$.

We look for the rank-2 vector bundles, which appears as extensions of type:
$$0\longrightarrow\mathcal{L}_1\longrightarrow\mathcal{E}\longrightarrow\mathcal{L}_2\longrightarrow 0,$$
where $\mathcal{L}_1=\mathcal{O}_{V^n}(-tE_p-tE_{2n+2})$ and $\mathcal{L}_2=\mathcal{O}_{V^n}(tE_p+tE_{2n+2})$.
For this extension, $H^1(V^n,{\mathcal{L}_2}^{\vee}\otimes\mathcal{L}_1)=H^1(V^n,\mathcal{O}_{V^n}(-2tE_p-2tE_{2n+2}))$ is  equal due to Serre duality to the group
$$H^{n-1}(V^n, \mathcal{O}_{V^n}(K_{V^n}+2tE_p+2tE_{2n+2}))^{\vee},$$
which is exactely $H^{n-1}(V^n, \mathcal{O}_{V^n}(D))^{\vee}$, for $D=K_{V^n}+2t(E_p+E_{2n+2})$.

From Proposition \ref{omologia}, $$H^{n-1}(V^n, \mathcal{O}_{V^n}(D))=\bigoplus_{\textbf{r}\in\mathbb{Z}^{2n+2}}H_{1}(Supp(\textbf{r}), \mathbb{C}),$$
where the direct sum is taken over all the representations $\textbf{r}\in\mathbb{Z}^{2n+2}$ of the line bundle $\mathcal{O}_{V^n}(\sum_{i=1}^{2n+2} r_i E_i)\cong\mathcal{O}_{V^n}(D)$, and $H_{1}$ is the first group of reduced homology of the simplicial complex $\emph{Supp}(\textbf{r})$.
We are going to provide a description of all the representations $D=\sum_{i=1}^{2n+2} r_i E_i$ of our divisor, using the next lemma:
\begin{lemma}\label{Chow} The Chow ring of $V^n$ is equal to $\mathbb{Z}[E_i]_{1\leq i \leq 2n+2}/I$, where $I$ is the ideal generated by the relations:

\bigskip
$\bullet  E_{i_1}\cdot...\cdot E_{i_k} \cdot E_{j_1+n+1} \cdot...\cdot E_{j_{k'}+n+1} = 0,$ for all $i_1,...,i_k,j_1,...,j_{k'} \in \{1,...,n+1\}$, with the condition $k > [n/2]$ or $k' > [n/2]$ or $\{i_1,...,i_k\} \cap \{j_1,...,j_{k'} \} \neq \emptyset$.

\bigskip

$\bullet \sum_{i=1}^{2n+2}  \langle m,v_i \rangle E_i = 0$, for all $m \in Hom(N_{V^n}, \mathbb{Z})$.

\end{lemma}

\begin{proof} The proof follows directly from the description of the Chow ring associated to a toric variety, (see for instance, Proposition page 106, \cite{Fu}), and from the definition of the fan of $V^n$.
\end{proof}
From Lemma $\ref{Chow}$, we see that we have the linear relations: 
$$E_p-E_{n+1}+E_{2n+2}-E_{p+n+1}=0,$$
which describe the equivalence of divisors.
Using this relations, we get all the representations of $D$ in the Chow's ring:
$$O_{V^n}(D)=\mathcal{O}_{V^n}(K_{V^n}+2t(E_p+E_{2n+2}))$$
$$\,\,\,\,\,\,\,\,\,\,\,\,\,\,\,\,\,\,\,=\mathcal{O}_{V^n}(K_{V^n}+(2t-q)(E_p+E_{2n+2})+q(E_{n+1}+E_{p+n+1})), \,\textrm{where}\,\,q\in\mathbb{Z}.$$
Since $K_{V^n}=-\sum_{k=1}^{2n+2}E_k$, it follows that the representations of the divisor $D$ take the form:
\begin{equation}\label{divizorul} D=-\sum_{k\neq p, n+1, p+n+1, 2n+2 }E_k+(2t-q-1)(E_p+E_{2n+2})+(q-1)(E_{n+1}+E_{p+n+1}).\end{equation}
We can now construct the simplicial complex attached to the positive coefficients of the above representation of the divisor $D$.

If $q<0$, then $2t-q-1>0$ and $q-1<0$. This implies that $D$ has only two coefficients bigger or equal with zero, the ones corresponding to the toric divisors $E_p$ and $E_{2n+2}$, associated to the rays generators  $v_p$ and $v_{2n+2}$.
So, the geometric realization encodes all the cones generated by these two ray generators. 
It results that the reduced complex associated to $Supp(\textbf{r})$ has:
\newline $C_2=0$, since it dont exist $3$-dimesional cones generated by $v_p$ and $v_{2n+2}$, we dont have a $2$-simplex.
\newline $C_1=\mathbb{C}$, because $\langle v_p, v_{2n+2}\rangle$ is the only $2$-dimensional cone generated by the rays generators $v_p$ and $v_{2n+2}$, and so $\langle v_p, v_{2n+2}\rangle$ is the only $1$-simplex from the simplicial complex.
\newline $C_0=\mathbb{C}^2$, both $\langle v_p\rangle$ and $\langle v_{2n+2}\rangle$ are $1$-dimensional cones in the fan of $V^n$, so our simplicial complex contains two $0$-simplices;
\newline Of course, $C_{-1}=\mathbb{C}$, since this group encodes the zero cone.
Therefore, the reduced chain complex associated to the simplicial complex $Supp(\textbf{r})$ corresponding to the divisor $D$, is given for $q<0$ by:
$$0\stackrel{\delta_2}{\longrightarrow}\mathbb{C}\stackrel{\delta_1}{\longrightarrow}\mathbb{C}^2\stackrel{\delta_0}{\longrightarrow}\mathbb{C}\stackrel{\delta_{-1}}{\longrightarrow} 0,$$
with the operators $\delta_i$ defined as follows:
\newline $\delta_{-1}$ the zero map;
\newline $\delta_0(\langle v_i\rangle)=1,$ $\forall$ $i=p$ or $2n+2$;
\newline $\delta_1(\langle v_{p}, v_{2n+2})=\langle v_p\rangle-\langle v_{2n+2}\rangle$.

Notice that $rank(\delta_1)=1$ implies $ker (\delta_1)=0$. On the other hand, $Im\delta_2=0$. 
Hence, we can conclude, that for $q<0$ the first reduced homology group $$H_1(Supp(\textbf{r}, \mathbb{C}))=\frac{ker\delta_1}{Im \delta_2}=0.$$

We will find the same result, for $q\geq 2t$. This time, the only positive coefficients that appear in the description of the divisor $D$ will be the ones corresponding to the toric divisors $E_{n+1}$ and $E_{p+n+1}$. As in the above case, we obtain the same reduced homology complex, and thus the same conclusion.

If $1\leq q\leq 2t-1$, than our divisor $D$ has four  coefficients bigger or equal to zero, the ones associated to the toric divisors:
$E_p, E_{n+1}, E_{n+p+1}$ and $E_{2n+2}$.
Thus, in this case the reduced complex of $Supp(\mathbf{r})$ is supported on the cones of the fan of $V^n$ generated by the rays generators $v_{p}$, $v_{n+1}$, $ v_{p+n+1}$ and $v_{2n+2}$. 

First of all, we remark  from the definition of the fan of $V^n$ that these four rays generators dont determinate $4$-dimensional cones or $3$-dimensional in our variety. Therefore, $C_3=C_2=0$ since we dont have $3$-simplices or $2$-simplices.

Following the same ideea as in the previous cases, we get:
\newline $C_1=\mathbb{C}^4$, since we have four $1$-simplices coresponding to the $2$-dimensional cones $\langle v_p, v_{n+1}\rangle$, $\langle v_p, v_{n+1}\rangle$, $\langle v_{n+1}, v_{p+n+1}\rangle,$ $\langle v_{p+n+1}, v_{2n+2}\rangle.$
\newline $C_0=\mathbb{C}^4$, because we have four $0$-simplices which encodes the four $1$-dimensional cones $\langle v_p\rangle $, $\langle v_{n+1}\rangle$, $\langle v_{p+n+1}\rangle$, $\langle v_{2n+2}\rangle$, and again
\newline $C_{-1}=\mathbb{C}$, for the zero cone.
\newline Hence, for $1\leq q\leq 2t-1$, the reduced chain complex of $Supp(\mathbf{r})$ is given by:
$$0\stackrel{\delta_2}{\longrightarrow}\mathbb{C}^4\stackrel{\delta_1}{\longrightarrow}\mathbb{C}^4\stackrel{\delta_0}{\longrightarrow}\mathbb{C}\stackrel{\delta_{-1}}{\longrightarrow} 0.$$

We remark that $\delta_1$ is defined as follows: $$\delta_1(\langle v_p, v_{n+1}\rangle)=\langle v_{n+1}\rangle-\langle v_p\rangle,$$
$$\delta_1(\langle v_p, v_{2n+2}\rangle)=\langle v_{2n+2}\rangle-\langle v_p\rangle,$$
$$\delta_1(\langle v_{n+1}, v_{p+n+1}\rangle)=\langle v_{p+n+1}\rangle-\langle v_{n+1}\rangle,$$
$$\delta_1(\langle v_{p+n+1}, v_{2n+2}\rangle)=\langle v_{2n+2}\rangle-\langle v_{p+n+1}\rangle.$$
Since the rank of  $\delta_1$ is the same as the rank of the matrix
$$\left(%
\begin{array}{cccc}
  -1 & -1 & 0 & 0 \\
   1 &  0 & -1& 0 \\
   0 &  0 &  1 & -1 \\
   0 &  1 &  0 & 1 \\
\end{array}%
\right),$$ we get $Im\delta_1=\mathbb{C}^3$. 
Therefore $ker\delta_1=\mathbb{C}$. Since $Im\delta_2=0$, we obtain $$H_1(Supp(\textbf{r}))=\frac{ker\delta_1}{Im \delta_2}=\mathbb{C}.$$
We remind that since $1\leq q\leq 2t-1$, we have $2t-1$ representations of the line bundle $\mathcal{O}_{V^n}(D)$ with $H_1(Supp(r))=\mathbb{C}$.
Hence,
$$H^1(V^n, \mathcal{O}_{V^n}(-2tE_{p}-2tE_{2n+2})\stackrel{Serre}{=}H^{n-1}(V^n,\mathcal{O}_{V^n}(D))\cong$$ $$\stackrel{Prop.\ref{omologia}}{\cong}\bigoplus_{q=1}^{2t-1}H_1(Supp(\textbf{r}))=\mathbb{C}^{2t-1}.$$
As a consequence of the previous discussion, we have the next result:
\begin{theorem}\label{neanulare}
Let $V^n$ be a $n>2$ dimensional Del Pezzo variety, $t\geq 1$ an integer and $E_i$ the toric divisors corresponding to the ray generators $v_i$ of $V^n$ , for $i\in\{1,...,2n+2\}$.

Then for every $p\in\{1,...,n\}$, we have the following
nonvanishing result:  
$$H^1(V^n,\mathcal{O}(-2tE_{p}-2tE_{2n+2}))={2t-1} .$$
\end{theorem}

\bigskip
\textit{ Proof of Theorem 4.1}
Now, let $D_1 = -tE_p-tE_{2n+2}$  with $t \geq 1$ and $D_2 = -D_1$. By  Theorem \ref{neanulare}, we
have a nontrivial extension $\mathcal{E}$ of $\mathcal{O}_{V^n}(D_1)$
by $\mathcal{O}_{V^n}(D_2)$. We will prove that for any other exact sequence:
$$0 \rightarrow \mathcal{O}_{V^n}(X_1) \rightarrow \mathcal{E}
\rightarrow \mathcal{O}_{V^n}(X_2) \rightarrow 0, $$ we have:
$$ X_1 \cong D_1 \,\, \textrm{and} \,\, X_2 \cong D_2,$$
or
$$ X_1 \cong -D_1 \,\, \textrm{and} \,\, X_2 \cong -D_2.$$ This implies that $\mathcal{E}$ is unsplit.

So, assume we have such an exact sequence. By Whitney's formula for
Chern classes we have:
$$c_1(\mathcal{E}) = X_1+X_2 = D_1 + D_2 = 0,$$ and
$$c_2(\mathcal{E}) = X_1 \cdot X_2 = D_1 \cdot D_2.$$
As a consequence, we get \begin{equation}\label{intersectia}X_2^2-D_2^2 =0.\end{equation} Using Lemma \ref{Chow} we analize the above relation.

Let suppose that $X_2=\sum_{i=1}^{2n+2} \alpha_iE_i$, $\alpha_i\in\mathbb{Z}$. 
By the relation (\ref{intersectia}), it follows
$$(\sum_{j=1}^{2n+2}{\alpha}_j^2E_j^2+2\sum_{l<j}{\alpha}_l{\alpha}_j E_l E_j)-(t^2E_p^2+t^2E_{2n+2}^2+2t^2E_pE_{2n+2})=0,$$
or  equivalent
$$\sum_{j=1, j\neq p,2n+2}^{2n+2}{\alpha}_j^2E_j^2+2\sum_{l<j, l\neq p,j\neq 2n+2}{\alpha}_l{\alpha}_j E_l E_j+$$
$$+({\alpha}_p^2-t^2)E_p^2+({\alpha}_{2n+2}^2-t^2)E_{2n+2}^2+2(\alpha_p\alpha_{2n+2}-t^2)E_pE_{2n+2}=0.$$
\c Now from the description of Chow ring (see Lemma \ref{Chow}), we get the relations:
\newline ${\alpha}_j^2=0$ for $j\neq p, 2n+2$;
\newline${\alpha}_p^2-t^2=0\Longleftrightarrow \alpha_p=\pm t$;
\newline${\alpha}_{2n+2}^2-t^2=0\Longleftrightarrow\alpha_{2n+2}=\pm t$;
\newline $\alpha_p\alpha_{2n+2}-t^2=0\Longrightarrow$ $\alpha_p$  and $\alpha_{2n+2}$ must have the same sign .(Beware that, $E_pE_{n+2}\neq 0$,  since there are cones in the fan of $V^n$ that contains $v_p$ and $v_{2n+2}$ in the same time.)
Therefore, the only options for $X_2$ are : $$X_2=D_2=tE_p+tE_{2n+2}$$ or $$X_2=-D_2=-tE_p-tE_{2n+2},$$
and thus the proof is completed.

\textbf{Acknowledgements} Part of this work was done while both authors were invited at the Universit\'a di Catania. It is a pleasure to thank Laura Costa and Rosa Maria Mir\'o-Roig for many inspiring discussions.
The second author also thanks to Roland Abuaf for insightful comments on the paper.

\begin{tabbing}
\hspace*{8cm} \= \hspace*{2cm} \= \hspace*{3cm} \= \hspace*{4cm}
\kill

\parbox{\textwidth}{Giulio Cotignoli}    \> \parbox{\textwidth}{Alexandru Sterian}\\
\parbox{\textwidth}{University of Rome 'La Sapienza'} \> \parbox{\textwidth}{University "Spiru Haret"}\\
\parbox{\textwidth}{Department of Mathematics}   \> \parbox{\textwidth}{Department of Mathematics}\\
\parbox{\textwidth}{P.le Aldo Moro, 5, 00185} \> \parbox{\textwidth}{Str. Ion Ghica, 13, 030045}\\
\parbox{\textwidth}{Rome, Italy}  \> \parbox{\textwidth}{Bucharest, Romania}\\
\parbox{\textwidth}{cotignoli@mat.uniroma1.it}   \> \parbox{\textwidth}{alexandru.sterian@gmail.com}\\

\end{tabbing}

\end{document}